\documentclass{article}
\usepackage{spconf,amsmath,amssymb,graphicx,algorithm,color}

\def \om {{\omega    }}

\def \qH {{\mathbb  H}}
\def \N  {{\mathbb  N}}
\def \R  {{\mathbb  R}}
\def \C  {{\mathbb  C}}

\def \cS {{\mathcal S}}
\def \cV {{\mathcal V}}
\def \cF {{\mathcal F}}
\def \cH {{\mathcal H}}
\def \cA {{\mathcal A}}
\def \sgn{\mathop{\rm sgn}}

\newcommand{\beqn}{\begin{equation}}
\newcommand{\eeqn}{\end{equation}}
\newcommand{\eref}[1]{(\ref{#1})}

\newtheorem{theorem}{Theorem}[section]

\newtheorem{defi}[theorem]{Definition}

\newtheorem{cor}[theorem]{Corollary}

\title{A unique polar representation of the hyperanalytic signal}
\name{Boqiang Huang$^{\star}$, Angela Kunoth$^{\star \dagger}$
      \thanks{B. Huang is grateful to the Alexander von Humboldt foundation.}}
\address{$^{\star}$ Institut f\"ur Mathematik,
         Universit\"at Paderborn, Germany\\
         E-mail: bhuang@math.upb.de \\
         $^{\dagger}$ Mathematisches Institut,
         Universit\"at zu K\"oln, Germany\\
         E-mail: kunoth@math.uni-koeln.de}

\begin{document}
\maketitle
\begin{abstract}
The hyperanalytic signal is the straight forward generalization
of the classical analytic signal. It is defined by a complexification
of two canonical complex signals, which can be considered as an
inverse operation of the Cayley-Dickson form of the quaternion.
Inspired by the polar form of an analytic signal where the real
instantaneous envelope and phase can be determined, this paper
presents a novel method to generate a polar representation of the
hyperanalytic signal, in which the continuously complex envelope and
phase can be uniquely defined. Comparing to other existing methods,
the proposed polar representation does not have sign ambiguity
between the envelope and the phase, which makes the definition of the
instantaneous complex frequency possible.
\end{abstract}
\section{Introduction}
\label{sec:intro}

In the signal processing community, the analytic signal is
computed as a complexification of a real-valued signal, generated
by the signal itself and its Hilbert transform, see, e.g.,
\cite{Picinbono1997}. This is a well-known model for signal
characterization. Based on it, the instantaneous amplitude or
envelope, the instantaneous phase, and, thus, the instantaneous
frequency of the given signal can be well identified. Then, with
the obtained time and instantaneous amplitude and frequency
quantities, the constructed time-frequency-amplitude (TFA)
spectrum may illustrate valuable information for data classification,
signal decomposition and many other applications.

Nowadays, in many applications, e.g., geophysical
\cite{Lilly2010, Rudi2010} or meteorological \cite{Rehman2010} data
analysis, signals simultaneously sampled from multiple sensors
may not be efficiently characterized based on the classical model.
Moreover, the subcomponents of non-stationary multivariate
signals by using an adaptive data decomposition method, e.g.,
\cite{Rehman2010, Huang2013}, also need a versatile TFA representation.
Therefore, it is necessary to develop a solid theory for multivariate
signal analysis.

A multivariate version based on quarternions, so-called monogenic
signals, with applications to images was developed in \cite{Jager2010}.
In \cite{Bihan2014}, the concept of hyperanalytic signal (H-signal) was
proposed to provide a hypercomplex representation of the given complex
or bivariate signal by using a one-sided quaternionic Fourier transform
\cite{Said2008}. Inspired by the Cayley-Dickson form, the generated
H-signal may be represented in a polar form where the corresponding
envelope and the phase are complex \cite{Sangwine2010}. However, there
is a sign ambiguity between the envelope and the phase, which results
in the fact that the definition of the instantaneous complex frequency
is unclear. \cite{Lilly2010} introduced another method for complex
signal characterization based on the modulated elliptical model
\cite{Schreier2008}.

This paper firstly reviews the quaternion computation and the
H-signal construction from a given complex signal. Then the polar
representation of the quaternion and the corresponding sign
ambiguity will be explained. After that we present a novel envelope
recovery algorithm based on a linear zero-crossing prediction,
which results in an unique polar form for the H-signal. Thus, the
instantaneous complex frequency can be naturally defined. We illustrate
the efficiency of the proposed method via a representative numerical
study, and close with a discussion and some final remarks.

\section{Quaternionic signals}
\label{sec:signal}

\subsection{Quaternion and its operations}
\label{ssec:operation}

A quaternion, introduced by Hamilton in 1843, is for the Cartesian
coordinate axes the set
$\qH := \{q_r + iq_i + jq_j + kq_k \; : \; q_r,q_i,q_j,q_k \in \R\}$
where $\{1,i,j,k\}$ customarily denotes the basis. Every element
$q \in \qH$ can be uniquely written in a linear combination of these
basis elements.

The real part of $q$, denoted as the scalar, is $\cS(q):=q_r$, while
the residual is denoted by $\cV(q) := q - \cS(q)$. $q$ is called a
pure quaternion when $q_r = 0$. Each basis element is considered as
the root of $-1$. Multiplications among them satisfy
$i^2 = j^2 = k^2 = ijk = -1$, which results in the potential
rules $ij = -ji = k, jk = -kj = i, ki = -ik = j$. Note that the
quaternion multiplication is not commutative, e.g., $qp \neq pq,$
for $q,p \in \qH$. The conjugate of $q$ is defined as
$\bar{q} = q_r - iq_i - jq_j - kq_k$, and, thus, the Euclidean
norm of $q$ is
$\|q\| : = \sqrt{q \bar{q}} = \sqrt{q_r^2+q_i^2+q_j^2+q_k^2}$. Then
the inverse of $q$ is given by
$q^{-1} := \frac{\bar{q}}{\|q\|^2}$. By applying the Cayley-Dickson form,
any quaternion can be represented as a pair of complex numbers, e.g.
$q = q_r + iq_i + jq_j + kq_k = (q_r + iq_i) + (q_j + iq_d)j =: z_1 + z_2 j$,
for $z_1, z_2 \in \C$.

\begin{defi}\label{Def_ExpQuat}
Given a quaternion $q \in \qH$, the exponential and the natural
logarithm of $q$ can be defined by
\begin{small}
\begin{eqnarray}
e^{q} &:=& \textstyle e^{\cS(q)} \left( \cos(\|\cV(q)\|) +
           \frac{\cV(q)}{\|\cV(q)\|}\sin(\|\cV(q)\|) \right),
           \label{Eq_ExpQuat}\\
\ln(q)&:=& \textstyle \ln(\|q\|) +
           \frac{\cV(q)}{\|\cV(q)\|}\arccos(\frac{\cS(q)}{\|q\|}).
           \label{Eq_LogQuat}
\end{eqnarray}
\end{small}
\end{defi}

\subsection{Hyperanalytic signal}
\label{ssec:transform}

For any complex signal $z(t) \in \C$, the quaternionic Fourier
transform (QFT) of $z(t)$ may have left, right and double-sided
versions, since the exponential kernel placed in different positions
leads to different results. In the present context,
the right QFT is appropriate
for the H-signal construction \cite{Bihan2014, Said2008}.

\begin{defi}\label{Def_RightQFT}
Given a complex signal $z(t) \in \C$, and a unit quaternion
$\mu \in \qH$, the right QFT of $z(t)$ with respect to (w.r.t)
the $\mu$-axis is defined by
\beqn\label{Eq_RightQFT_1}
\hat{z}_\mu(\om)  = \cF_{\mu}^{q}[z(t)]
                 := \textstyle \int_{\R} z(t) e^{-\mu \om t} dt
\eeqn
\end{defi}

If we replace the $\mu$-axis with a canonical $j$-axis, the
$\cF_j^q$ of a real signal $a(t)$ can be considered as the Fourier
transform of $a(t)$. Thus, the QFT can be sped up by the fast Fourier
transform (FFT) as follows.

\begin{cor}\label{Cor_RightQFT}
Given a complex signal
$z(t) = z_r(t) + iz_i(t)$, $t, z_r(t), z_i(t) \in \R$, and the
quaternionic $j$-axis, the right QFT of $z(t)$ can be expressed in
terms of the Fourier transform
\beqn\label{Eq_RightQFT_2}
\hat{z}_j(\om) = \cF_j^{q}[z(t)] = \cF_j[z_r(t)] + i \cF_j[z_i(t)]
\eeqn
\end{cor}

Based on the QFT, we can modify the Hilbert transform (HT) in
the Hamilton space. In the time domain, the output of the HT of a
complex signal should still be a complex signal that is orthogonal to
the input. Then this pair of signals can be combined in a
Cayley-Dickson form to generate a quaternionic signal (Q-signal). In
the frequency domain, the frequency of such generated Q-signal should
be physically meaningful, i.e., nonnegative.

\begin{defi}\label{Def_HT}
Given a complex signal $z(t) \in \C$ and a unit quaternion
$\mu \in \qH$, the quaternionic Hilbert transform (QHT) of $z(t)$
w.r.t the $\mu$-axis is defined by
\beqn\label{Eq_HT_1}
\cH_{\mu}^{q}[z(t)] := {\cF_{\mu}^{q}}^{-1}
          \left[-\mu \sgn(\om) \cF_{\mu}^{q}[z(t)] \right].
\eeqn
\end{defi}
Here, ${\cF_{\mu}^{q}}^{-1}$ means the inverse QFT. The QHT can
also be defined in the time domain by
$\cH_{\mu}^{q}[z(t)] := \text{PV} (z(t)\ast\frac{1}{\pi t})$, where
PV denotes the Cauchy principal value and $\ast$ represents the
convolution. Replacing again the $\mu$-axis with the $j$-axis,
\eref{Eq_HT_1} can be further simplified.

\begin{cor}\label{Cor_RightQFT_new}
Given a complex signal
$z(t) = z_r(t) + iz_i(t)$, $t, z_r(t), z_i(t) \in \R$, and the
quaternionic $j$-axis, the QHT of $z(t)$ can be expressed in
terms of the HT
\beqn\label{Eq_HT_2}
\cH_j^{q}[z(t)] = \cH [z_r(t)] + i \cH [z_i(t)].
\eeqn
\end{cor}

Similar to the analytic signal model, we can construct the H-signal
for any given complex signal, which is indeed a subset of the
Q-signal.

\begin{defi}\label{Def_Hsignal}
Given a complex signal $z(t) \in \C$, the hyperanalytic signal
is defined by
\beqn\label{Eq_HAsignal}
s(t) := z(t) + o(t)j = z(t) + \cH_j^{q}[z(t)] j,
\eeqn
where $o(t)$ is the QHT of $z(t)$ w.r.t the $j$-axis.
\end{defi}

\section{Hyperanalytical signal model}
\label{sec:HAmodel}

\subsection{Sign ambiguity in the polar form}
\label{sec:signproblem}

Suppose the quaternion $q = q_r + iq_i + jq_j + kq_k$ is given,
$q_r,q_i,q_j,q_k \in \R$, and its polar representation is in the
form of $q := A e^{Bj}$, where $A := a + ib, B := c + id,$ and
$a,b,c,d \in \R$. As $Bj$ is a pure quaternion, according to
\eref{Eq_ExpQuat}, the exponential of $Bj$ can be expressed as
\beqn\label{Eq_Bj}
\begin{small}
\begin{aligned}
e^{Bj} &:= \alpha + j\beta + k\gamma \\
       &:= \textstyle
           \cos(\|B\|) + j \frac{c}{\|B\|} \sin(\|B\|)
                       + k \frac{d}{\|B\|}\sin(\|B\|),
\end{aligned}
\end{small}
\eeqn
where $\|B\|=\sqrt{c^2+d^2}$ \cite{Sangwine2010}. Then, we arrive
at the equations
\beqn\label{Eq_qEquality}
\begin{aligned}
q &= q_r + i q_i + j q_j + k q_k := Ae^{Bj}  \\
  &= a \alpha + i b \alpha +
     j (a \beta - b \gamma) + k (a \gamma + b \beta).
\end{aligned}
\eeqn
Since the complex envelope $A$ can be expressed in polar form by
$A := \|A\|e^{i \phi_A} = \|q\|e^{i \phi_A}$, we can determine that
the axis of the known complex component $q_r + i q_i$ equals to the
axis of $a \alpha + i b \alpha$. In other words, with an axis
operator defined as $\cA(a+ib):=\frac{a+ib}{\|a+ib\|}$, we
have following relationship
\beqn\label{Eq_AxisEq}
e^{i \phi_A} = \cA(a+ib) = \textstyle \frac{\cA(q_r + i q_i)}{\sgn(\alpha)},
\eeqn
where $\sgn(\cdot)$ is the signum function. Obviously, this leads
to an ambiguity in sign between the complex envelope $A$ and the
phase $B$ since the $\sgn(\alpha)$ is unknown for computing
$e^{i \phi_A}$.

\subsection{Complex envelope recovery}
\label{sec:envelope}

To simplify the polar representation of the H-signal, denoted by
$q(t)=A(t)e^{B(t)j}, t \in [0, T]$, we assume a unit signal to be
processed in this section, i.e., $\|q(t)\|=1$. Then the unwanted
$\sgn(\alpha)$ in \eref{Eq_AxisEq} can be removed by taking the
modulus of the real and imaginary components on both sides,
\beqn\label{Eq_AxisAbs}
\begin{aligned}
|\cos(\phi_A(t))| &= |a(t)| = |\tilde{q}_r(t)|, \\
|\sin(\phi_A(t))| &= |b(t)| = |\tilde{q}_i(t)|,
\end{aligned}
\eeqn
where $\cA(q_r(t) + i q_i(t)) := \tilde{q}_r(t) + i\tilde{q}_i(t)$.

Recall that we are considering a continuous hyperanalytical signal
model, in which the complex envelope $A(t)$ and, thus, the real
phase $\phi_A(t)$ should be continuous. Therefore, if the initial range
of the phase was limited, i.e., $\phi_A(0) \in [0, \frac{\pi}{2}]$,
in view of the continuity, the recovered envelopes $\tilde{a}(t)$
and $\tilde{b}(t)$ could be determined independently, as
$|\tilde{q}_r(t)|$ and $|\tilde{q}_i(t)|$ are already known. In
detail, since the sign changing of the envelope $a(t)$ or $b(t)$ only
occurs at zero-crossing (ZC) position which is nothing but the local
minimum of the modulus $|\tilde{q}_r(t)|$ or $|\tilde{q}_i(t)|$, in
principle, we can recover the envelope by retrieving the sign of
every half-period (HP) of the modulus signal from beginning to the
end, where the HP is defined as the interval between every two nearest
local minima of the modulus signal.

\begin{figure}[t]
\begin{minipage}[b]{.48\linewidth}
  \centering
  \centerline{\includegraphics[width=4.7cm]{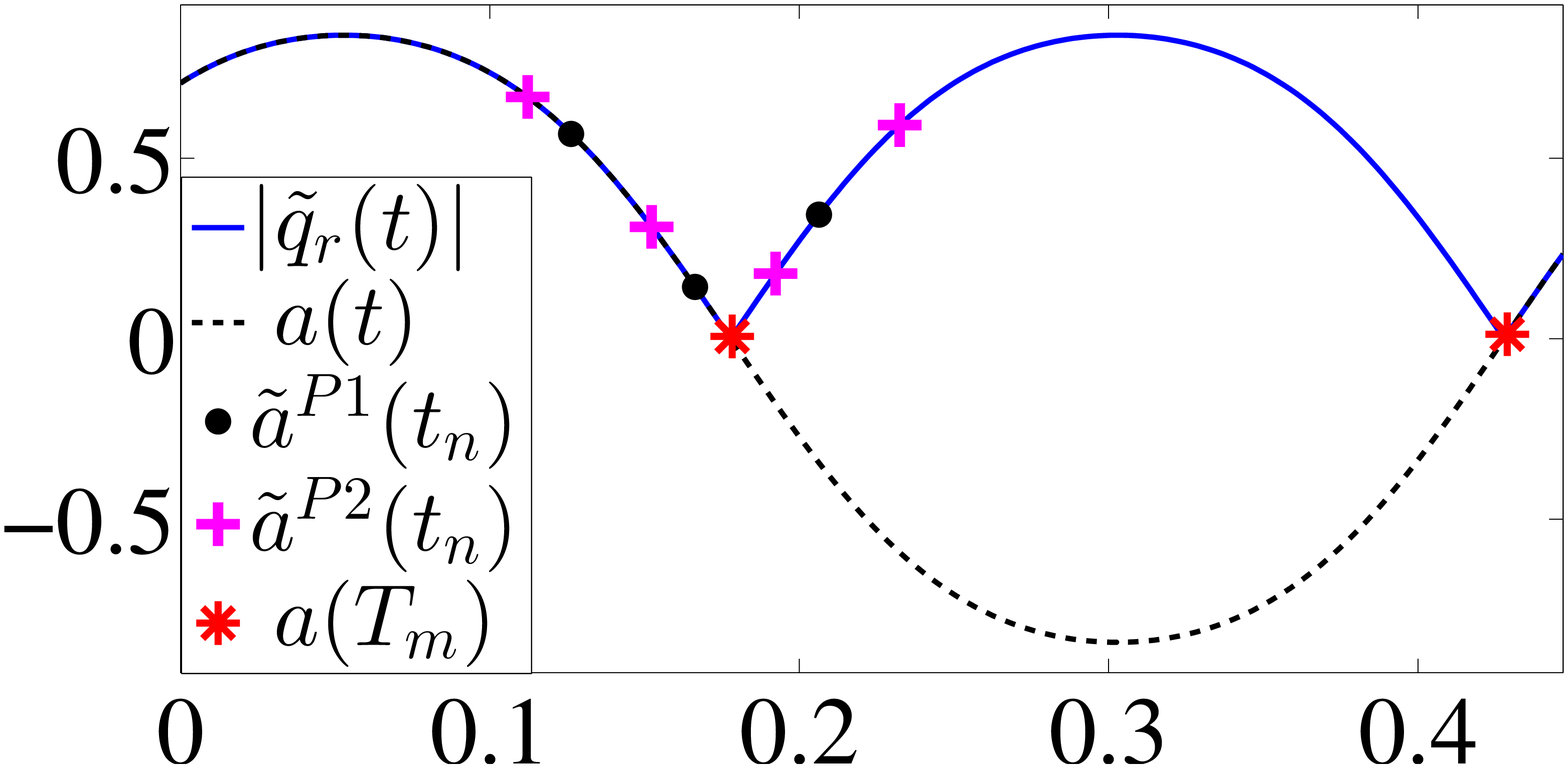}}
  \centerline{\small (a) Class I with cases $P1$, $P2$}\medskip
\end{minipage}
\hfill
\begin{minipage}[b]{.48\linewidth}
  \centering
  \centerline{\includegraphics[width=4.7cm]{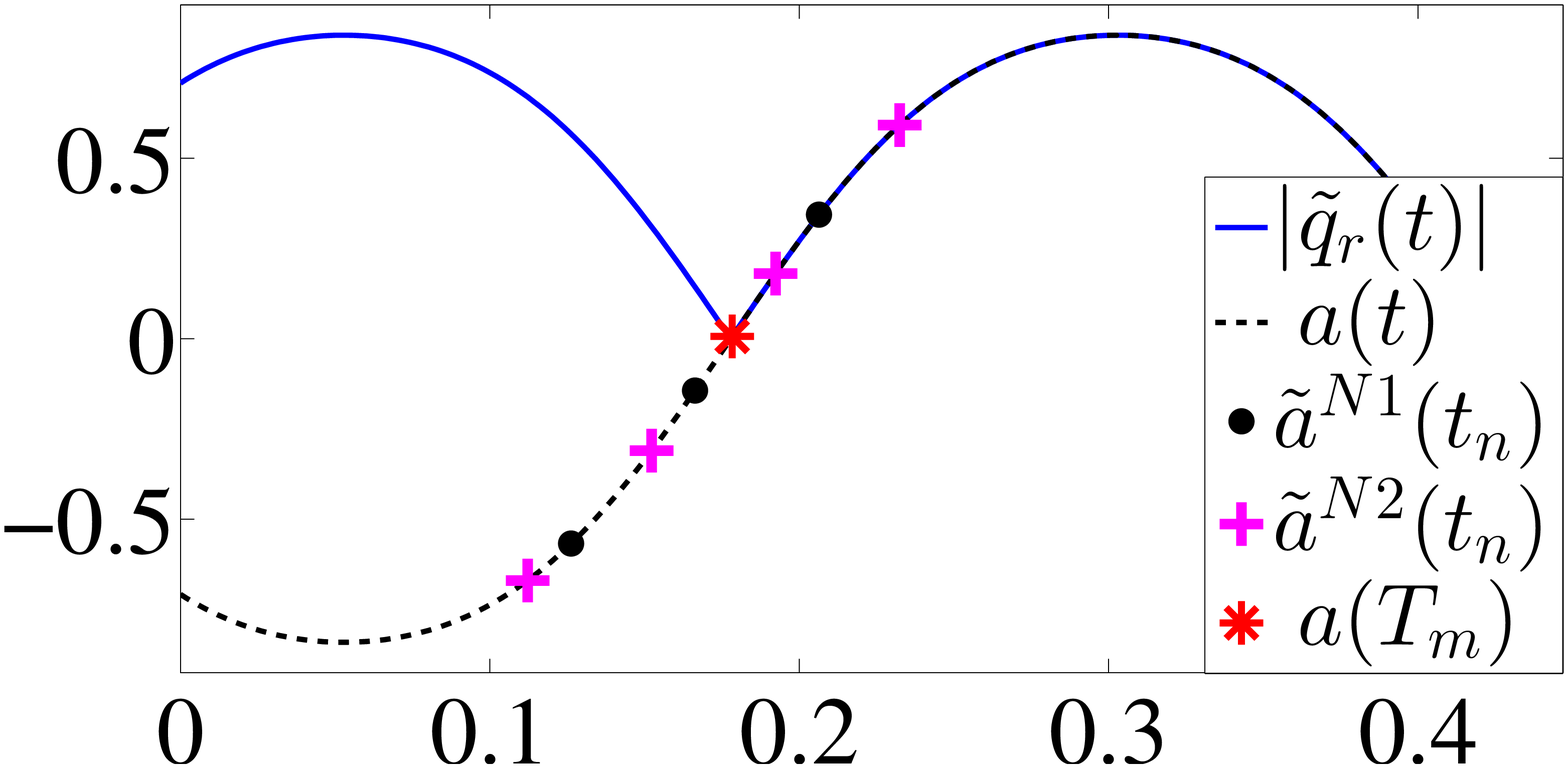}}
  \centerline{\small (b) Class I with cases $N1$, $N2$}\medskip
\end{minipage}
\begin{minipage}[b]{.48\linewidth}
  \centering
  \centerline{\includegraphics[width=4.7cm]{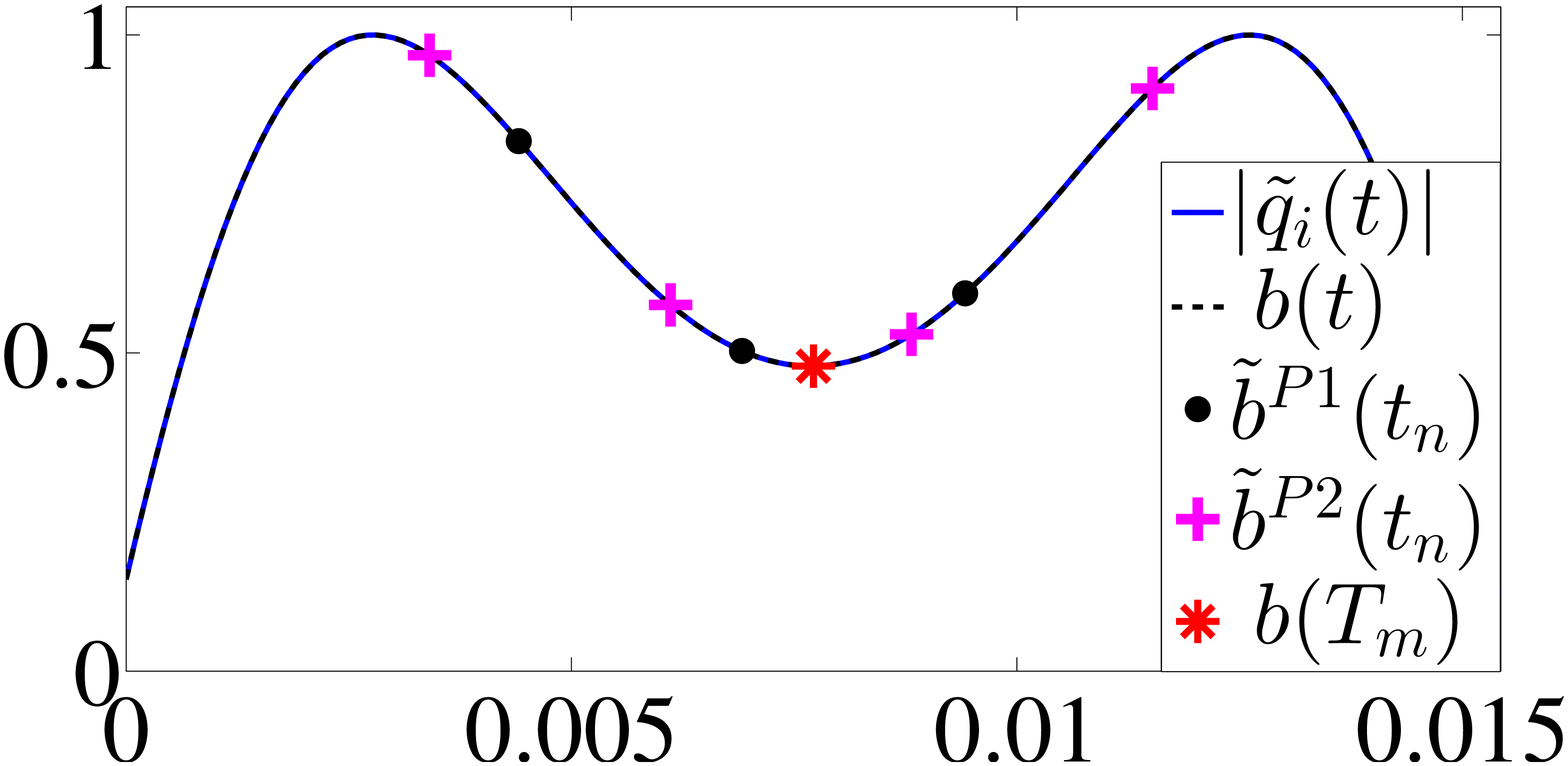}}
  \centerline{\small (c) Class II with cases $P1$, $P2$}\medskip
\end{minipage}
\hfill
\begin{minipage}[b]{.48\linewidth}
  \centering
  \centerline{\includegraphics[width=4.7cm]{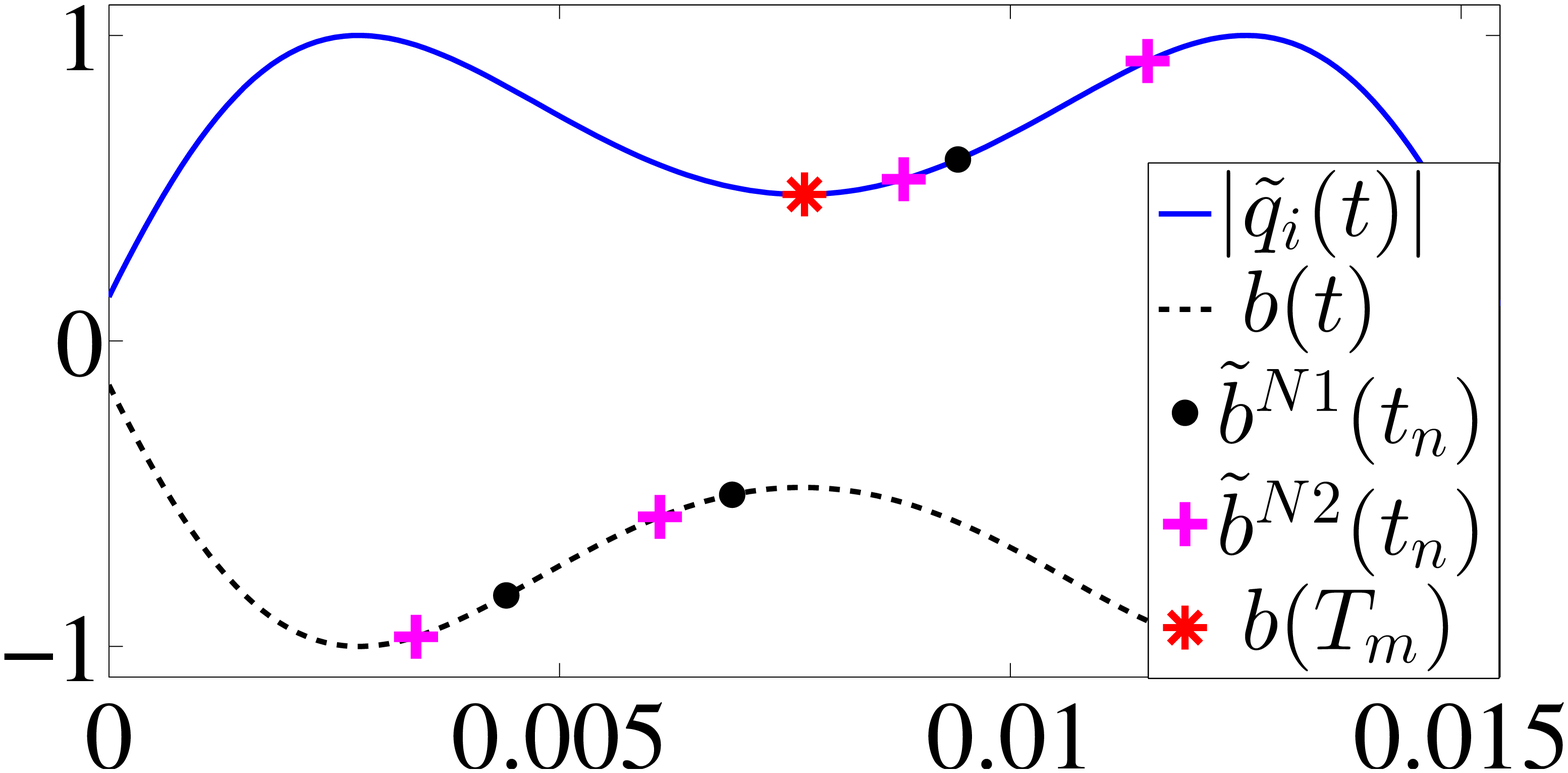}}
  \centerline{\small (d) Class II with cases $N1$, $N2$}\medskip
\end{minipage}
\caption{\small Eight possibilities for the sign recovery of
the envelope $A(t)=a(t)+ib(t)$ based on the modulus components of
the axis $\cA(q_r(t) + i q_i(t))$. $P$ (positive) and $N$ (negative)
denote the sign of the former half-period of the recovered component
$\tilde{a}(t)$ and $\tilde{b}(t)$. Black dots and magenta pluses
denote sampling cases 1 and 2, and red stars imply the ideal local
minima (Color online).}
\label{fig:EnvRec}
\end{figure}

However, the local minimum of the modulus signal may be positive but
not the ZC because we do not require that the phase $\phi_A(t)$ is
monotonically non-decreasing. Therefore, we need to classify all cases
into two classes: class I denotes the case the local minimum is the
ZC, while class II implies a positive local minimum. In addition,
for discrete data, the accuracy of the local minimum position is
affected by the sampling rate, which means that the current local
minimum may be the last point of the former HP (case 1), or the first
point of the following one (case 2). Therefore, if we ignore special
cases for stationary points, in total, there will be eight possibilities
which may occur around the local minimum. For instance, in class I, we
have to consider the former HP is positive (case $P$) or negative
(case $N$), each of which contains another two sampling cases.
Fig.~\ref{fig:EnvRec} gives a comprehensive illustration of all
possibilities.

To distinguish these different possibilities, we employ a linear ZC
prediction method based on every two successive samples around the
local minimum. Take the case $P1$ in Fig.~\ref{fig:EnvRec} (a) as
an example. Black dots are denoted by
$\tilde{a}(t_{n-1}), \tilde{a}(t_n)$ and $\tilde{a}(t_{n+1})$, among
which $t_n$ corresponds to the local minimum. The predicted ZC
position is
\beqn\label{Eq_ZCPredict}
T_{\tilde{a}}^n := t_{n+1} - \tilde{a}(t_{n+1})
\textstyle\frac{t_{n+1} - t_n}{\tilde{a}(t_{n+1}) - \tilde{a}(t_n)}.
\eeqn
Similarly, $T_{\tilde{b}}^{n}$ can be calculated for the imaginary
component. Then, with the information of the sign of the samples
$\tilde{a}(t_{n-1})$ and $\tilde{a}(t_{n+1})$ and the estimates
$T_{\tilde{a}}^n$ and $T_{\tilde{a}}^{n-1}$, we can determine to
which case the current local minimum belongs. Also considering
the black dots in Fig.~\ref{fig:EnvRec} (a), we can firstly determine
that the former HP is positive (case $P$) as
$\sgn(\tilde{a}(t_{n-1})) = \sgn(\tilde{a}(t_{n+1}))$.
Secondly, we can determine the class I as the estimate
$T_{\tilde{a}}^{n-1}$ is valid ($t_n \leq T_{\tilde{a}}^{n-1} \leq t_{n+1}$),
and simultaneously the sampling case 1 as the estimate
$T_{\tilde{a}}^n$ is invalid ($T_{\tilde{a}}^n < t_{n-1}$). Therefore,
we can keep the sign of $\tilde{a}(t_n)$ and then change the sign of
the following NP from $t_{n+1}$ to the former point of the next local
minimum.

Since the case determination contains many IF--ELSE conditions, we
only present a simplified envelope recovery algorithm as follows.
The accuracy of the ZC prediction is guaranteed as the sampling
frequency is high enough, otherwise the instantaneous frequency
cannot be correctly estimated because of the violation of the
sampling theorem.
\begin{algorithm}
\renewcommand{\thealgorithm}{:}
\caption{Complex Envelope Recovery Algorithm}
\small
\label{Alg_EnvRec}
\begin{itemize}
\item[1.] Initialize the recovered components
    $\tilde{a}(t) = |\tilde{q}_r(t)|$ ($\tilde{b}(t) = |\tilde{q}_i(t)|$),
    and the range of $\phi_A(0)$, e.g., $\phi_A(0) \in [0, \frac{\pi}{2}]$;
\item[2.] Detect the local minimum of $\tilde{a}(t)$ ($\tilde{b}(t)$), and
    retrieve the sign of the first HP based on the range of the phase;
\item[3.] Recover the envelope by retrieving the sign of the flowing HP
    based on the sign of the former HP and the determined case at the local
    minimum $\tilde{a}(t_n)$ ($\tilde{b}(t_n)$);
\item[4.] Output the complex envelope $A(t) = \tilde{a}(t)+i\tilde{b}(t)$.
\end{itemize}
\end{algorithm}

\subsection{Unique polar representation}
\label{sec:Polar}

Once the complex envelope is recovered, the quaternionic carrier
$e^{B(t)j}$ can be computed by
\beqn\label{Eq_Carrier}
e^{B(t)j} := \alpha(t) + j\beta(t) + k\gamma(t)
           = \textstyle\frac{\bar{A}(t)q(t)}{\|q(t)\|^2},
\eeqn
and the complex phase $B(t)$ can be derived based on \eref{Eq_LogQuat}
\beqn\label{Eq_PhaseB}
B(t) := c(t) + id(t) = \cA(\beta(t) + i\gamma(t)) \arccos(\alpha(t)).
\eeqn
From \eref{Eq_Bj}, we know that there is still a sign ambiguity between
$\sin(\|B(t)\|)$ and $c(t)$ or $d(t)$. However, since both $\sin(\cdot)$
and $\cos(\cdot)$ functions are periodic, we have that
$\sin(\|B(t)\|) = \sin(\|B(t)\| \pm 2m \pi)$, and
$\cos(\|B(t)\|) = \cos(\|B(t)\| \pm 2m \pi)$, $m \in \N$. Thus, it
is reasonable to assume that $c(t), d(t)$ are non-negative and
monotonically non-decreasing, and the initial phase $c(0), d(0)$ should
satisfy $\|B(0)\|:=\sqrt{c(0)^2+d(0)^2} \in [0, 2 \pi)$. Finally, we can
uniquely retrieve the phases $c(t)$ and $d(t)$ based on the unwrapped
$\arccos(\alpha(t))$ in \eref{Eq_PhaseB}. The reason to unwrap the phase
$\arccos(\alpha(t))$ but not the ones $\breve{c}(t)$ and $\breve{d}(t)$
(which are directly computed in \eref{Eq_PhaseB} without unwrapping) is
because only the phase $\arccos(\alpha(t))$ has a fixed period $2\pi$.
Then the retrieved phases $\tilde{c}(t)$ and $\tilde{d}(t)$ can be
considered as approximations of the ideal ones that are monotonically
non-decreasing. Thus, we have proved the following result.

\begin{theorem}\label{Theorem_HSignal}
Given a complex signal $z(t) \in \C$, the hyperanalytic signal can be
constructed by $s(t) := z(t) + \cH_j^{q}[z(t)] j$, $s(t) \in \qH$, which
has an unique polar form $s(t) = A(t)e^{B(t)j}$, $A(t), B(t) \in \C$, if
$(A(t), B(t))$ is the canonical complex pair where
$A(t):=\|s(t)\|e^{i \phi_A(t)}$, $\phi_A(0) \in [0, \frac{\pi}{2}]$,
and $B(t):=c(t)+id(t)$, $c(t),d(t) \geq 0$, $\|B(0)\| \in [0, 2\pi)$.
\end{theorem}
Bearing in mind that the instantaneous frequency should be nonnegative,
therefore, each component of the unwrapped complex phase $B(t)$ should be
monotonically non-decreasing, and the unwrapped phase $\phi_A(t)$ is
the same if the complex envelope $A(t)$ is an analytic signal.

\begin{defi}\label{Def_IF}
Given a complex signal $z(t) \in \C$ and let the polar form of its
hyperanalytic signal be
$s(t) := A(t)e^{B(t)j} := \|s(t)\|e^{i \phi_A(t)}e^{(c(t)+id(t))j}$,
$\phi_A(t),c(t),d(t) \geq 0$. The instantaneous complex frequency of
$z(t)$ is defined by
\beqn\label{Eq_ICF}
f_B(t) := f_{B_r}(t) + i f_{B_i}(t)
=\textstyle \frac{1}{2 \pi} (\frac{d(c(t))}{dt} + i\frac{d(d(t))}{dt}),
\eeqn
and the instantaneous frequency of the complex envelope $A(t)$ is defined by
$f_A(t) := \frac{d(\phi_A(t))}{2 \pi dt}$ if $A(t)$ is an analytic signal.
\end{defi}

\section{Numerical study}
\label{sec:NumricalStudy}

To illustrate the efficiency of the proposed method, we design a
representative hyperanalytical signal model
\beqn\label{Eq_Numerical}
s(t) := e^{-t} e^{7\sin(2 \pi t)i}
        e^{(40 \pi t + i(20 \pi t + 4\cos(2 \pi t)))j},
\eeqn
for $t \in [0, 0.4]$, from which one can define
$A(t) := a(t) + i b(t) = e^{-t} e^{7\sin(2 \pi t)i}$, and
$B(t) := c(t) + i d(t) = 40 \pi t + i(20 \pi t + 4\cos(2 \pi t))$.
Thus we can obtain the given complex signal $z(t)=z_r(t)+iz_i(t)$ based on
\eref{Eq_HAsignal} and \eref{Eq_qEquality}, and determine the instantaneous
complex frequency by
$f_B(t) := f_{B_r}(t) + i f_{B_i}(t) = 20 + i (10 - 4 \sin(2 \pi t))$.
Here, $f_A(t)$ is not well defined since it can be negative.

\begin{figure}[t]
\begin{minipage}[b]{.48\linewidth}
  \centering
  \centerline{\includegraphics[width=4.4cm]{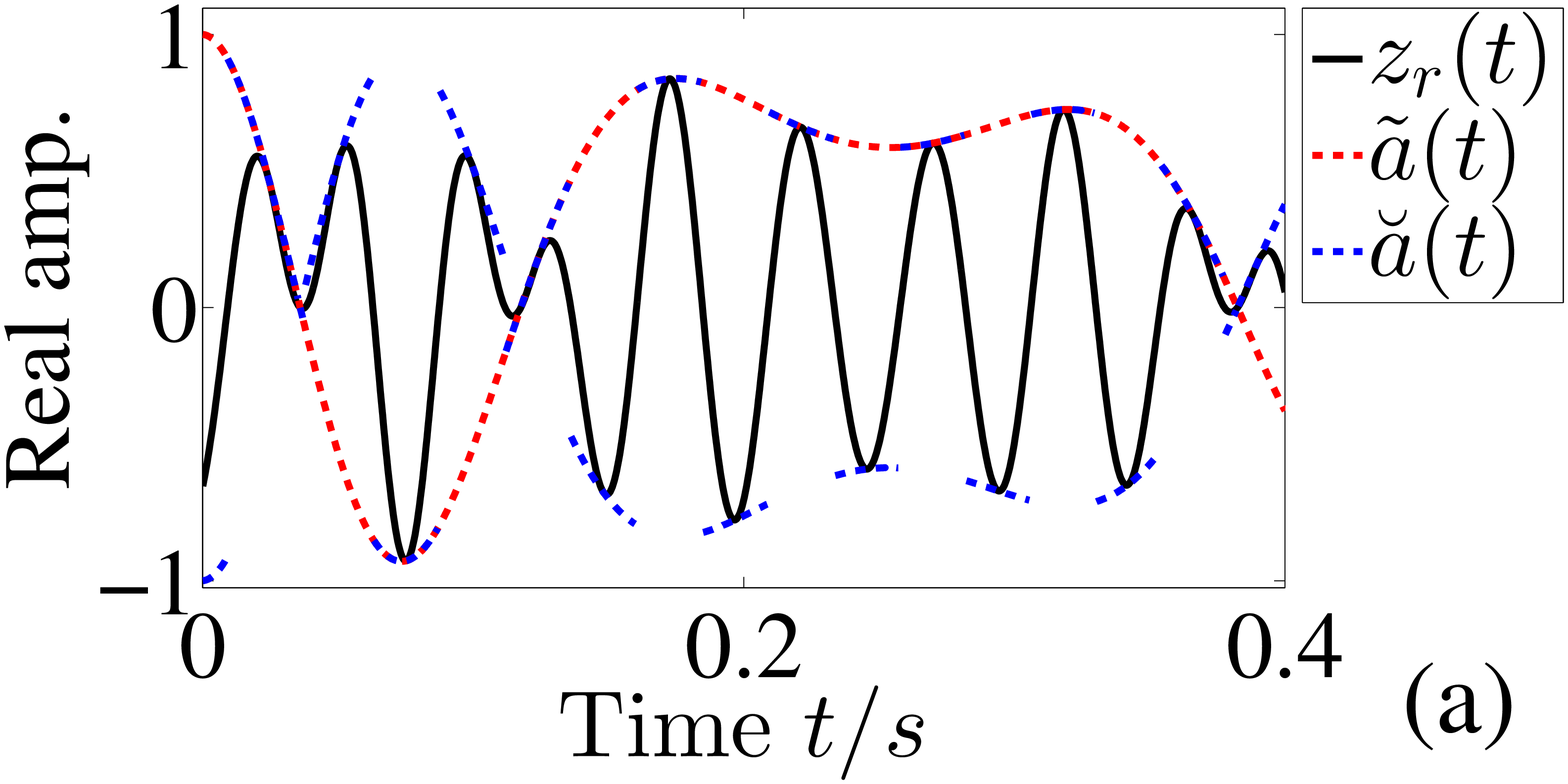}}
\end{minipage}
\hfill
\begin{minipage}[b]{.48\linewidth}
  \centering
  \centerline{\includegraphics[width=4.4cm]{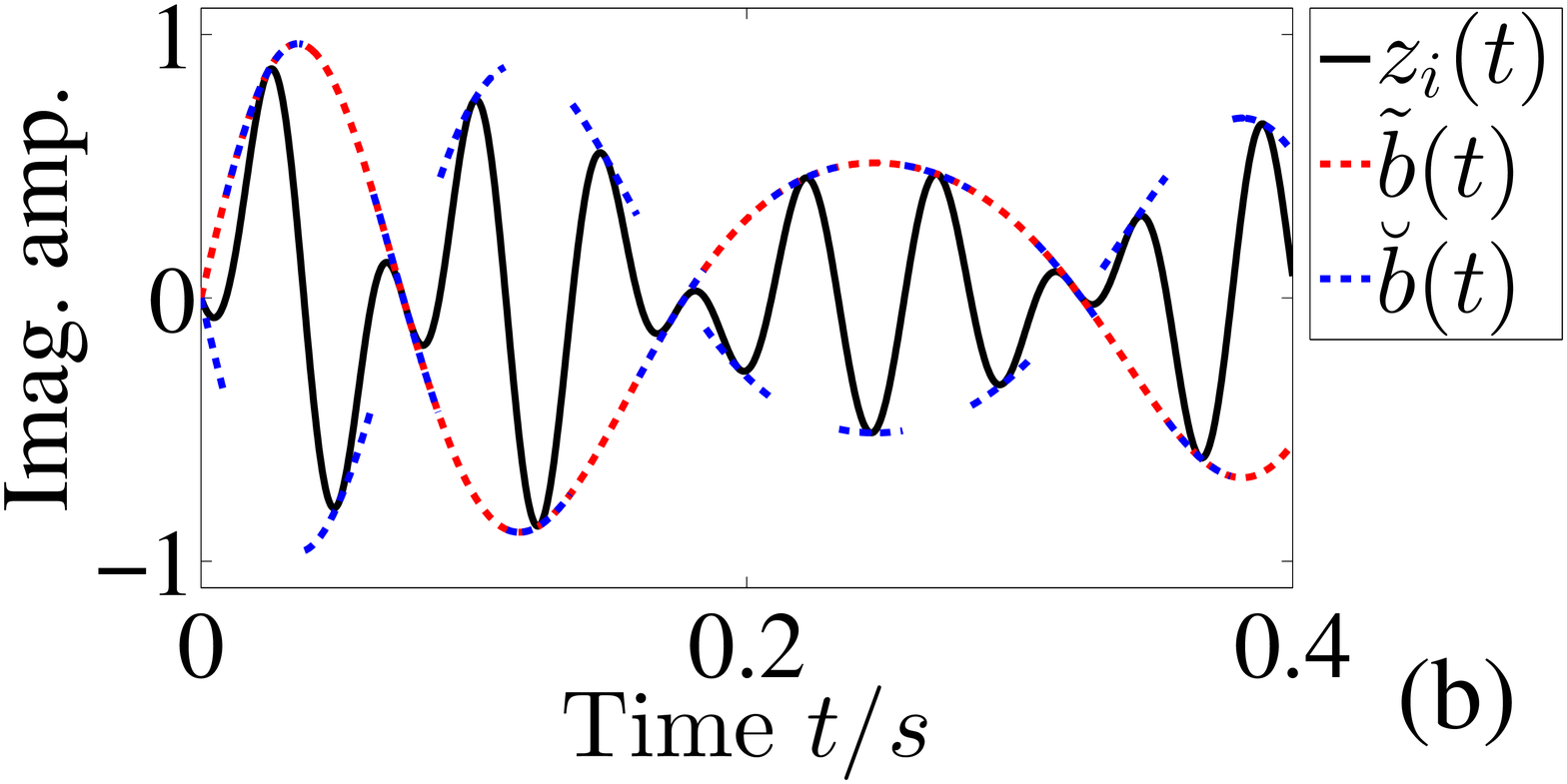}}
\end{minipage}
\begin{minipage}[b]{.48\linewidth}
  \centering
  \centerline{\includegraphics[width=4.4cm]{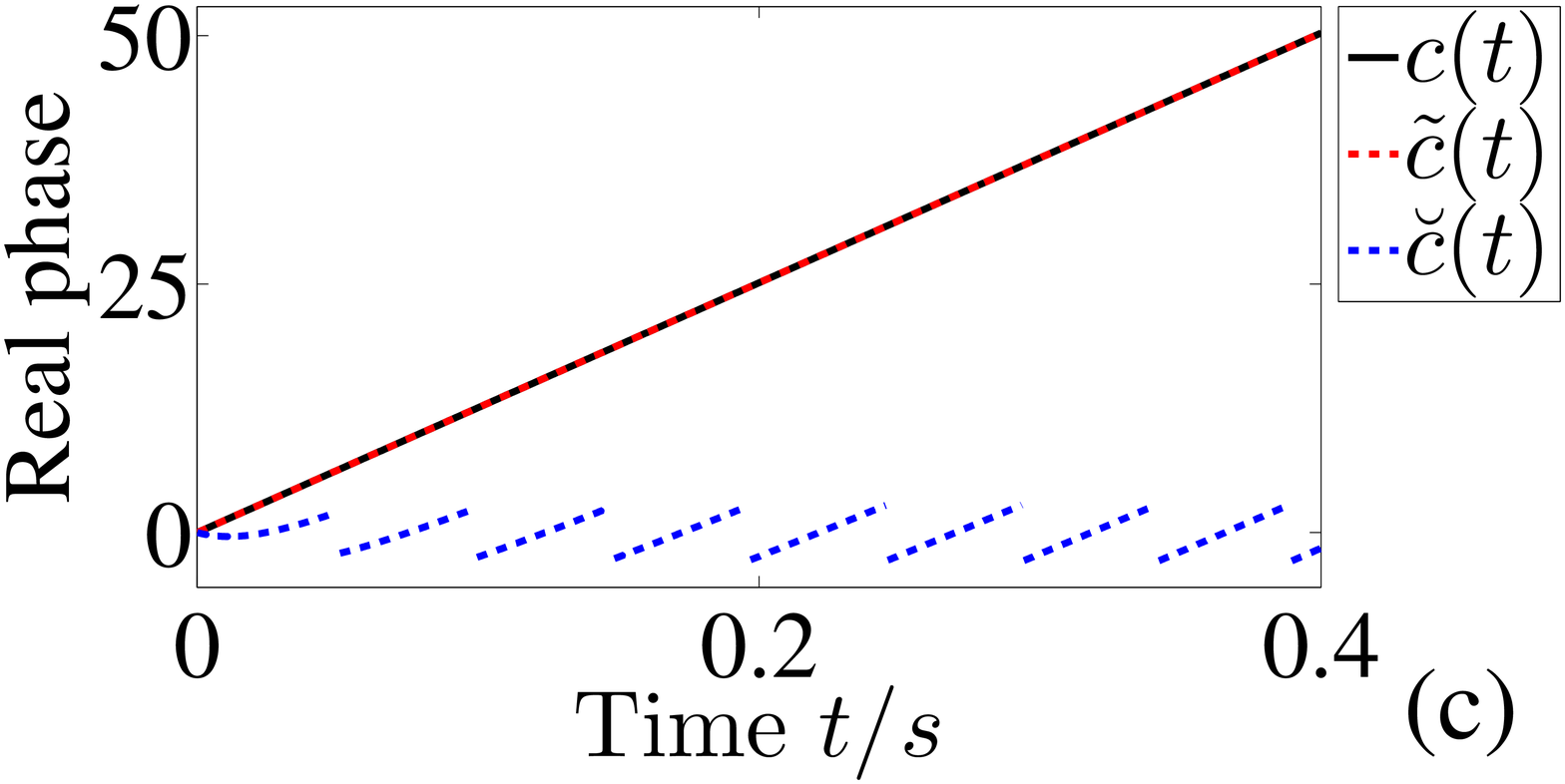}}
  \medskip
\end{minipage}
\hfill
\begin{minipage}[b]{.48\linewidth}
  \centering
  \centerline{\includegraphics[width=4.4cm]{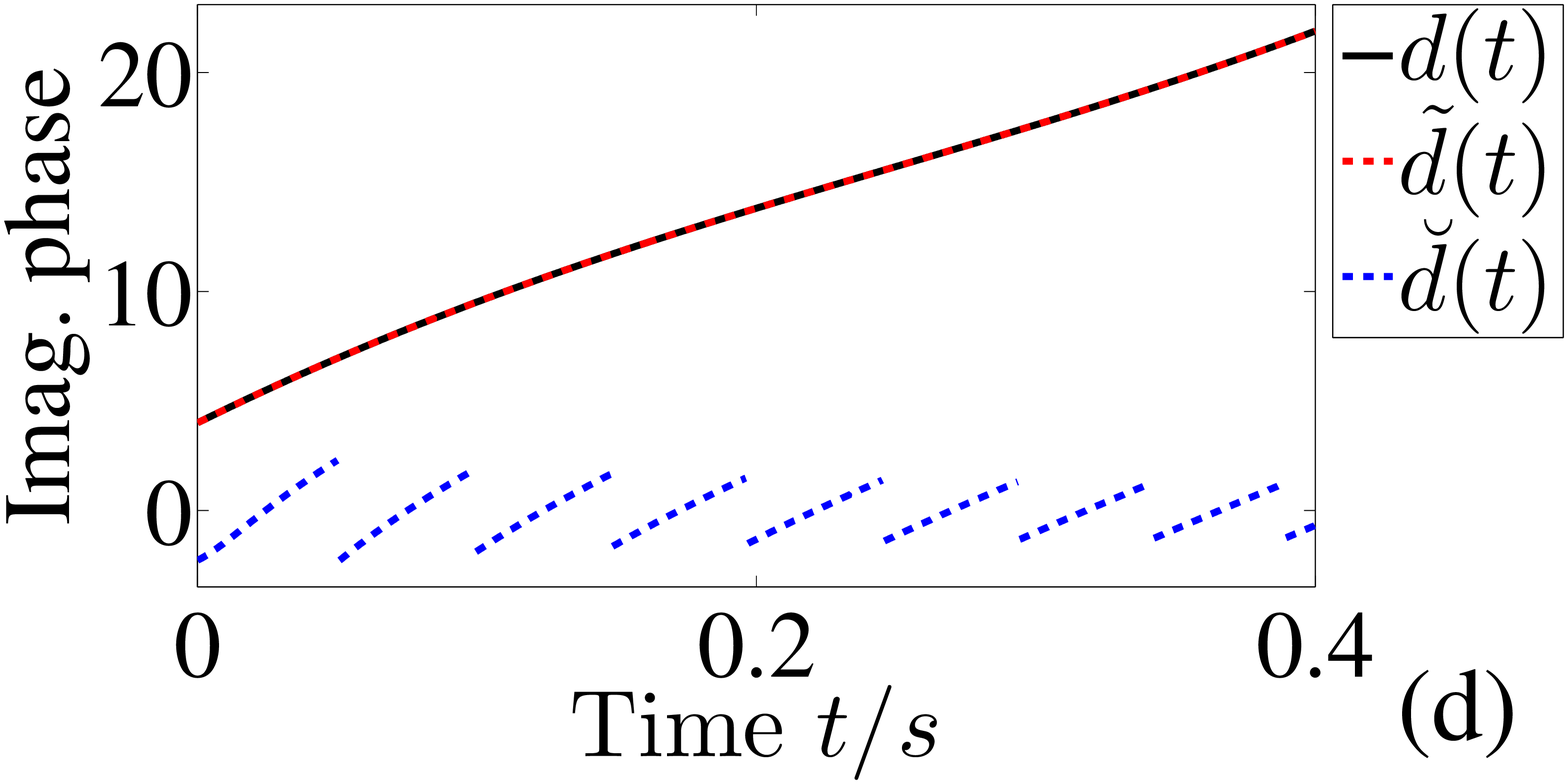}}
  \medskip
\end{minipage}
\begin{minipage}[b]{.48\linewidth}
  \centering
  \centerline{\includegraphics[width=4.4cm]{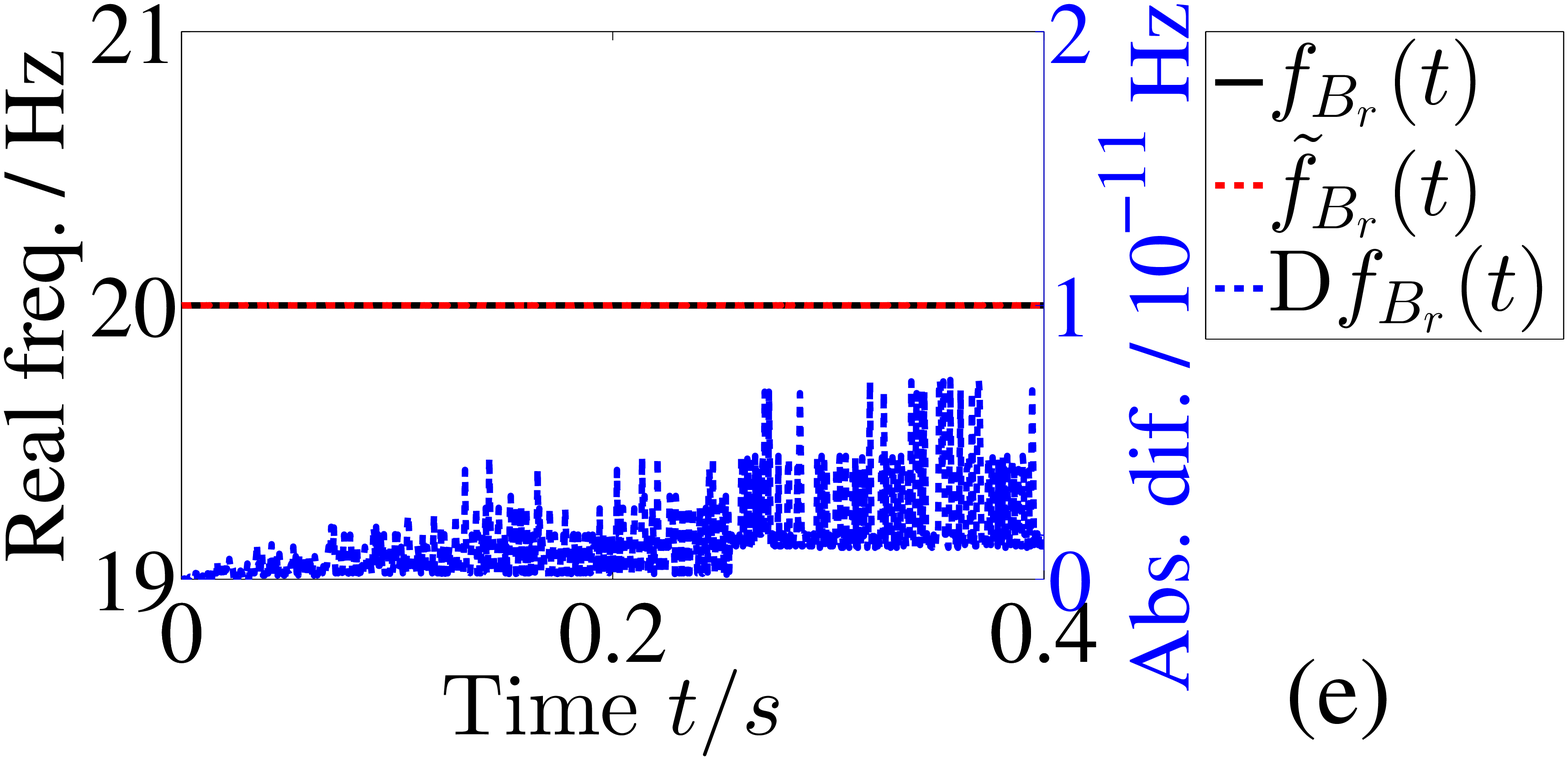}}
  \medskip
\end{minipage}
\hfill
\begin{minipage}[b]{.48\linewidth}
  \centering
  \centerline{\includegraphics[width=4.4cm]{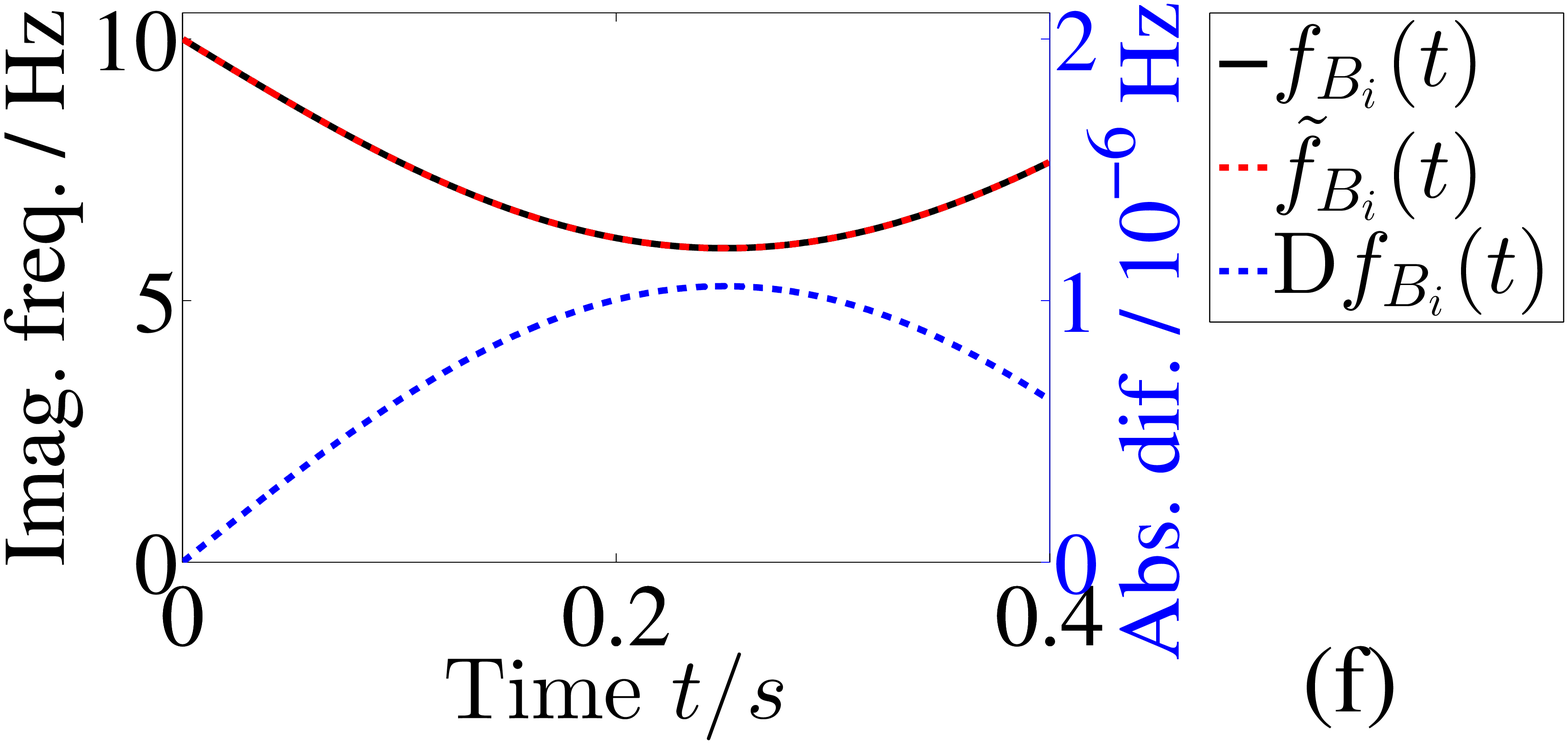}}
  \medskip
\end{minipage}
\caption{\small Numerical results for the H-signal model.
Top: real and imaginary (R\&I) parts of the given complex signal
$z(t)$, the recovered envelope $\tilde{A}(t)$, and the one $\breve{A}(t)$
containing the sign ambiguity; Middle: R\&I parts of the ideal complex
phase $B(t)$, the recovered one with
$(\tilde{B}(t) := \tilde{c}(t) + i \tilde{d}(t))$ and without
$(\breve{B}(t) := \breve{c}(t) + i \breve{d}(t))$ unwrapping;
Bottom: R\&I parts of the ideal instantaneous complex frequency $f_B(t)$,
the estimated one $\tilde{f}_B(t)$ based on $\tilde{B}(t)$, and the
absolute difference between each of them,
$\text{D}f_{B_r}(t):=|f_{B_r}(t)-\tilde{f}_{B_r}(t)|$,
$\text{D}f_{B_i}(t):=|f_{B_i}(t)-\tilde{f}_{B_i}(t)|$ (Color online).}
\label{fig:Numeric}
\end{figure}

Fig.~\ref{fig:Numeric} illustrates all respective results for the given
signal $z(t)$. In sub-figures (a) and (b),
$\breve{A}(t):=\breve{a}(t)+i \breve{b}(t) = \| q(t) \| \cA(q_r(t)+i q_i(t))$
is the reconstructed envelope based on \eref{Eq_qEquality}, which contains
the sign ambiguity. Obviously, the recovered
$\tilde{A}(t):=\tilde{a}(t)+i \tilde{b}(t)$
using the proposed method coincides strongly with the ideal complex envelope.
Sub-figures (c) and (d) imply the importance of the phase unwrapping,
while sub-figures (e) and (f) show the efficiency of the estimation of
the instantaneous complex frequency. Since the real component of
$f_B(t)$ is a constant, the absolute difference between it and the
estimated one $\tilde{f}_{B_r}$ is around machine accuracy. However,
the absolute difference between $f_{B_i}$ and $\tilde{f}_{B_i}$ is larger
since $f_{B_i}$ is nonlinear and thus the corresponding estimation
accuracy is corrupted by the discrete derivative computation at different
time positions.

\section{Conclusion}
\label{sec:Conclude}

We presented an efficient method for the unique polar representation of
the hyperanalytic signal that is constructed from any given complex
signal with continuous real and imaginary components. Based on this
H-signal model, we can obtain a canonical pair of continuously
instantaneous complex envelope and phase, in which the phase consists of
monotonically non-decreasing sub-components that leads to a natural
definition of the instantaneous complex frequency. Moreover, the
instantaneous real frequency of the complex envelope can also be
well-defined if the envelope is an analytical signal.

The developed H-signal model implies an interesting extension of the
multivariate signal characterization to arbitrary space dimensions, which
may have potential applications in such fields where the
time-frequency-amplitude information is representative for multivariate
signal analysis.

\bibliographystyle{IEEEbib}
\bibliography{Mystr}

\begin{thebibliography}{10}

\bibitem{Picinbono1997}
B.~Picinbono,
\newblock ``On instantaneous amplitude and phase of signals,''
\newblock {\em IEEE Trans. Sig. Proc.}, vol. 45, pp. 552--560, Mar. 1997.

\bibitem{Lilly2010}
J.M. Lilly and S.C. Olhede,
\newblock ``Bivariate instantaneous frequency and bandwidth,''
\newblock {\em IEEE Trans. Sig. Proc.}, vol. 58, pp. 591--603, Feb. 2010.

\bibitem{Rudi2010}
J.~Rudi, R.~Pabel, G.~Jager, R.~Koch, A.~Kunoth, and H.~Bogena,
\newblock ``Multiscale analysis of hydrologic time series data using the
  hilbert-huang-transform (hht),''
\newblock {\em Vadose Zone Journal}, vol. 9, pp. 925--942, Nov. 2010.

\bibitem{Rehman2010}
N.~Rehman and D.P. Mandic,
\newblock ``Empirical mode decomposition for trivariate signals,''
\newblock {\em IEEE Trans. Sig. Proc.}, vol. 58, pp. 1059--1068, Mar. 2010.

\bibitem{Huang2013}
B.~Huang and A.~Kunoth,
\newblock ``An optimization based empirical mode decomposition scheme,''
\newblock {\em J. Comput. Appl. Math.}, vol. 240, pp. 174--183, Mar. 2013.

\bibitem{Jager2010}
G.~Jager, R.~Koch, A.~Kunoth, and R.~Pabel,
\newblock ``Fast empirical mode decompositions of multivariate data based on
  adaptive spline-wavelets and a generalization of the hilbert-huang-transform
  (hht) to arbitrary space dimensions,''
\newblock {\em Adv. Adaptive Data Anal.}, vol. 2, pp. 337--358, July 2010.

\bibitem{Bihan2014}
N.Le Bihan, S.J. Sangwine, and T.A. Ell,
\newblock ``Instantaneous frequency and amplitude of orthocomplex modulated
  signals based on quaternion fourier transform,''
\newblock {\em Signal Proc.}, vol. 94, pp. 308--318, Jan. 2014.

\bibitem{Said2008}
S.~Said, N.Le Bihan, and S.J. Sangwine,
\newblock ``Fast complexified quaternion fourier transform,''
\newblock {\em IEEE Trans. Sig. Proc.}, vol. 56, pp. 1522--1531, Apr. 2008.

\bibitem{Sangwine2010}
S.J. Sangwine and N.Le Bihan,
\newblock ``Quaternion polar representation with a complex modulus and complex
  argument inspired by the cayley-dickson form,''
\newblock {\em Adv. in Appl. Clifford Algebras}, vol. 20, pp. 111--120, 2010.

\bibitem{Schreier2008}
P.J. Schreier,
\newblock ``Polarization ellipse analysis of nonstationary random signals,''
\newblock {\em IEEE Trans. Sig. Proc.}, vol. 56, pp. 4330--4339, Sep. 2008.

\end{thebibliography}

\end{document}